\input vanilla.sty
\pageheight {8.4in}
\pagewidth {6.9in}
\TagsOnRight
\pageno=1
\font\id=cmbx12  
\font\vb=cmb10
\def\n{\noindent}

\def\m{\medpagebreak}\def\br{\bigpagebreak}

\def\I{\varphi}
\def\f{\frac}

\def\E{\text{\bf E}}

\def\be{\beta}

\def\p{\partial}

\def\br{\bigpagebreak}

\define\om{\omega}

\def\be{\beta}

\def\p{\partial}\define\pa{\prime}

\font\bx=cmbx10

\define\cor#1{\m\n{\vb Corollary #1}. }
\define\Om{\Omega}
\define\ov{\overline}

\def\SS{\text{\bf S}}

\define\R{\text{\bf I\!\bf R}}

\define\ff#1{\f{\p^2f}{\p{#1}^2}}
\define\\E{\text{\bf I\!\bf E}}
\define\ipp{\int_0^{2\pi}}\define\ip{\int_0^{\f\pi2}}

\def\CB{\Cal B}

\font\dubi=cmssi12
\define\BBB{{\bold B}}

\vskip20mm

\n{\id  RECONSTRUCTION OF CENTRALLY SYMMETRICAL CONVEX BODIES}

\n{\id  BY PROJECTION CURVATURE RADII}

\vskip8mm

\n{\bf R. H. Aramyan}

\font\bx=cmbx10

\vskip8mm

\n{\lineskip11.381102pt\bx The article considers the problem
of existence and uniqueness of centrally symmetrical convex body
for which the projection curvature radius function
coincides with a given flag function.
A necessary and sufficient condition is found that ensures a positive
answer.
An algorithm for construction the body in question is proposed.}

\vskip18mm

\baselineskip22.762pt

\n{\bf\S 1. INTRODUCTION}

\n Let $F(\om)$ be a function
defined on the sphere $\SS^2$. The existence and
uniqueness of convex body $\BBB\subset\R^3$
for which the mean curvature
radius at a point on $\p\BBB$ with normal direction $\om$
coincides with given $F(\om)$ was posed
by Christoffel (see [4]). Let $R_1(\om)$ and $R_2(\om)$ be
the principle curvature radii of the surface of the body at the point
with normal $\om\in\SS^2$.
The Christoffel problem asked about the existence of $\BBB$
for which $R_1(\om)+R_2(\om)=F(\om)$.
The corresponding problem for
Gauss curvature $R_1(\om)R_2(\om)=F(\om)$ was posed and solved by Minkovski.
Blashke reduced the Christoffel problem to a
partial
differential equation of second order for support function (see [4]).
Aleksandrov and Pogorelov generelized these problems, and proved the
existence and
uniqueness of a convex body for which
$$G(R_1(\om),R_2(\om))=F(\om),\tag1.1$$
for a class of symmetric functions G (see [4], [6]).

\n In this paper we consider a similar problem posed
for the projection curvature radii
of
centrally symmetrical convex bodies (see [2]).
By $\CB_o$ we denote
the class of convex bodies $\BBB\subset\R^3$ that
have a center of symmetry at the origin $O\in\R^3$.
We use the notation:

\n $\SS^2$ -- the unit sphere in $\R^3$ (the space of
spatial directions),

\n $\SS_\om\subset\SS^2$ --
the great circle with pole at $\om\in\SS^2$,

\n $e(\om,\psi)$ -- the plane containing the origin
of $\R^3$ and the directions $\om\in\SS^2$ and
$\psi\in\SS_\om$,

\n $\BBB(\om,\psi)$ -- projection of $\BBB\in\CB_o$
onto $e(\om,\psi)$,

\n $R(\om,\psi)$ -- curvature radius of $\p\BBB(\om,
\psi)$ at the point whose outer normal direction is
$\om$.

\n Let $F(\om,\psi)$ be a symmetric
function defined on the space of "flags" $\{(\om,\psi):\,
\om\in\SS^2,\,\psi\in\SS_\om\}$ (see
[1]).
We pose the problem
of existence and
uniqueness of a convex body for which
$$R(\om,\psi)=F(\om,\psi),\tag1.2$$
and find a necessary and sufficient
condition on
$F(\om,\psi)$ that ensures a positive answer.
Note, that uniquness (up to parallel shifts) follows from
the classical uniqueness result on
Christoffel problem.

\n Now we describe the main result.
Let $F(\om,\psi)$ be a
function defined on the $\{(\om,\psi):\,
\om\in\SS^2,\,\psi\in\SS_\om\}$.
We say that the function $F(\om,\psi)$ is symmetric if
$F(\om,\psi)=F(\om,\pi+\psi)=F(-\om,\psi)$.
We define
$$\ov F(\om)=\f1\pi\int_0^{2\pi}F(\om,\psi)\,d\psi.\tag1.3$$
Below we use the usual spherical coordinates
$\nu,\phi$ on $\SS^2$
based on a choice $\om$ for the North Pole and a choice of
a reference point $\phi=0$ on the equator $\SS_\om$ (so points
$(0,\phi)$
lie on the equator $\SS_\om$). The point with coordinates
$\nu,\phi$ in that coordinates system
we will denote by $(\nu,\phi)_\om$

\proclaim{Theorem 1.1} A nonnegative symmetric
continuously differentable
function $F(\om,\psi)$ defined on the space

\n $\{(\om,\psi):\,
\om\in\SS^2,\,\psi\in\SS_\om\}$
represents the
projection curvature radius of some convex body
if and only if
$$F(\om,\psi)=\f{\ov F(\om)}2-\f1{2\pi}\int_0^{\f\pi2}
\int_0^{2\pi}\f{\ov F((\nu,\phi)_\om)\,
\cos2(\phi-\psi)}{\cos\nu}
\,d\phi\,d\nu\tag1.4$$
for all $\om\in\SS^2$ and all $\psi\in\SS_\om$ 
(integration order is important).
\endproclaim

\n Our proof of Theorem 1.1
suggests an algorithm
of construction of the convex body $\BBB$ for given $F(\om,\psi)$.
We will need the following facts from the convexity theory.

\n{\bf\S 2. PRELIMINARIES}

\m\n It is well known (see [8]) that
{\dubi the support function of every sufficiently smooth
$\BBB\in\CB_o$ has the unique representation
$$
H(\xi)=\int_{\SS^{2}}\mid(\xi,\Om)\mid h(\Om)\,d\Om,\tag2.1
$$
where $d\Om$ is the usual area measure on
$\SS^{2}$, $h(\Om)$ is an even continuous function
(not necessarily nonnegative) called the {\dubi
generating density} of $\BBB$.}

\n Below we will use the following result
by N. F. Lindquist, see [8].

\n {\dubi An even continuous function $h(\Om)$ defined on $\SS^2$
is the generating density of a body
$\BBB\in\CB_o$
if and only if
$$
\int_{\SS_\om}\cos^2(\widehat{\psi,\I})\,
h_\om(\I)\,d\I\ge0\tag2.2
$$
for all $\om\in\SS^2$ and $\psi\in\SS_\om$.}

\n Remarkably, the integral (2.2) has a clear geometrical
interpretation.

\n In [2] it was proved that {\dubi for any
sufficiently smooth $\BBB\in\CB_o$ and $\psi
\in\SS_\om$ the projection
curvature radius can be calculated as
$$
R(\om,\psi)=2\int_{\SS_\om}\cos^2(\widehat{\psi,\I})\,
h_\om(\I)\,d\I,\tag2.3
$$
where $(\widehat{\psi,\I})$ is the angle
between $\psi$ and $\I$, while $h_\om(\I)$
is the restriction of $h(\Om)$ of $\BBB$
on $\SS_\om$.}

\n In Blaschke's book [5] one can find the following representation
of the generating density.
{\dubi For any $\Om\in\SS^2$
$$
2\pi h(\Om)=\f1{4\pi}\int_{\SS^\Om}(R_1+R_2)\,d\tau-
\f1{4\pi}\ip\ipp (R_1((u,\tau)_\Om)+R_2((u,\tau)_\Om)^\pa_u
\f{1-\sin u}{\sin u}\,du\,d\tau,\tag2.4
$$
where $R_i((u,\tau)_\Om),\,i=1,2$  are the principal
curvature radii of $\BBB$ at the point with normal
$(u,\tau)_\Om$ (which has the spherical coordinates $u,\tau$ with
respect $\Om$)}

\n{\bf\S 3. PROOF OF THEOREM 1.1 AND CONSTRUCTION OF CONVEX BODY}

\m\n{\dubi Proof.} Necessity: let
$R(\om,\psi)$ be
the projection curvature radius of some convex body
$\BBB\in\CB_o$.
We have to
prove that $R(\om,\psi)$ satisfies the condition (1.4).
It follows from [5] that
$$\ov R(\om)=\f1\pi\int_0^{2\pi}R(\om,\psi)\,d\psi=
R_1(\om)+R_2(\om)\tag3.1$$
and
using Fubini Theorem one can rewrite the expression
(2.4) in the form
$$8\pi^2 h(\Om)=\lim_{a\to 0}\left[\f1{\sin a}
\int_0^{2\pi}\ov R((a,\tau)_\Om)\,d\tau-
\int_a^{\f\pi2}\int_0^{2\pi}\ov R((u,\tau)_\Om)
\f{\cos u}{\sin^2u}\,d\tau\,du\right],\tag3.2$$
where $(u,\tau)_\Om$ -is the point on $\SS^2$
with usual spheric coordinates
$u,\tau$ with respect to $\Om$.
Taking some $\psi\in\SS_\om$ for the reference point
on $\SS_\om$
and substituting the new expression for the
generating density
from (3.2) into (2.3) we get
$$R(\om,\psi)=2\int_{\SS_\om}\cos^2\I\,h_\om(\I)\,d\I=$$
$$=\lim_{a\to 0}\f1{4\pi^2}\int_0^{2\pi}
\cos^2\I\left[\f1{\sin a}
\int_0^{2\pi}\ov R((a,\tau)_\Om)\,d\tau-
\int_a^{\f\pi2}\int_0^{2\pi}\ov R((u,\tau)_\Om)
\f{\cos u}{\sin^2u}\,d\tau\,du\right]\,d\I,\tag3.3$$
where $\Om=(0,\I)_\om$.

\n For $u\in(0,\f\pi2)$ and $\I\in\SS_\om$
we denote by $S(u,\I)$
the circle
with center at $\I\in\SS_\om$ and of radius $\sin u$.
The meridian passing through the point $\om$ and $(0,\I)_\om$
divides $S(u,\I)$ into two halfcircles.
When $\I$ changes on $\SS_\om$, each
halfcircle of $S(u,\I)$ subtend
the set
$[0,2\pi)\times[-u,u]\subset\SS^2$. So on the latter set we have
two parametrizations
$(\nu,\phi)$ and $(\tau,\I)$.
One can prove the Jacobian relation
$$d\tau\,d\I=\f{\cos\nu}{\sqrt{\cos^2u-\sin^2\nu}}
d\nu\,d\phi.\tag3.4$$
Writing (3.3) as a sum of terms corresponding to the
halfcircles of $S(u,\I)$ and applying a change of variables (3.4),
we obtain
$$R(\om,\psi)=
\lim_{a\to 0}\f1{4\pi^2}\left[\f1{\sin a}\int_0^{2\pi}
\int_{-\f\pi2+a}^{\f\pi2-a}\ov R((\nu,\phi)_\om)
[\cos^2(\phi-\be)+\cos^2(\phi+\be)]
\f{\cos\nu}{\sqrt{\cos^2a-\sin^2\nu}}d\nu\,d\phi-\right.
$$
$$
-\left.\int_0^{2\pi}
\int_{-\f\pi2+a}^{\f\pi2-a}\ov R((\nu,\phi)_\om)
\left[\int_a^{\f\pi2-\nu}[\cos^2(\phi-\be)+
\cos^2(\phi+\be)]
\f{\cos\nu}{\sqrt{\cos^2u-\sin^2\nu}}
\f{\cos u}{\sin^2u} du\right]d\nu\,d\phi\right],\tag3.5
$$
where $\be=\phi-\I$. We have
$\cos\be=\f{\sin u}{\cos\nu}$
and $\sin\be=\f{\sqrt{\cos^2\nu-\sin^2u}}{\cos\nu}$.
After some
standard integral calculation (see [7])
and using symmetry we find
$$R(\om,\psi)=
\lim_{a\to 0}\f1{\pi^2}\left[\int_0^{2\pi}
\int_0^{\f\pi2-a}\ov R((\nu,\phi)_\om)
\cos^2\phi
\f{\sin a}{\cos\nu\sqrt{\cos^2a-\sin^2\nu}}d\nu\,d\phi-\right.
$$
$$
-\left.\int_0^{2\pi}
\int_0^{\f\pi2-a}\ov R((\nu,\phi)_\om)
\f{\cos2\phi}{\cos\nu}
\left[\f\pi2-\arcsin{\f{\sin a}{\cos\nu}}
\right]d\nu\,d\phi\right].\tag3.6
$$
Calculating the limit
requires decomposition of the integrals in the powers $a$.
The negative powers annihilate and we get
(1.4).

\n Sufficiently: let $F(\om,\psi)$ be a
nonnegative symmetric
continuous differentable function
satisfing the condition (1.4).
We consider $\ov F(\om)$ (see (1.3)) and
by means of (2.4) construct the function $f(\Om)$
defined on $\SS^2$:
$$f(\Om)=\f1{8\pi^2}\int_{\SS^\Om}\ov F((0,\tau)_\Om)\,d\tau-
\f1{8\pi^2}\ip\ipp \ov F((u,\tau)_\Om)^\pa_u
\f{1-\sin u}{\sin u}\,du\,d\tau.\tag3.7$$
According to N. F. Lindquist, $f(\Om)$ has to satisfy the condition
(2.2) to be the generating density of a body $\BBB$.
Substituting (3.7) into (2.2) and applying the same procedure
as above, by (1.4)
we obtaine $F(\om,\psi)$ which is nonnegative
by assumption.
Hence $h(\Om)$ is the generating function of $\BBB$
and according (2.3), $F(\om,\psi)$ is the projection
curvature radius of $\BBB$.

\n The author hopes to present soon the results on
similar problem that do not depend on the assumption
of central symmetry of the convex body $\BBB$.

\n I would like to express my gratitude to Professor
R. V. Ambartzumian for helpful discussions.

\baselineskip11.381102pt

\br\n{\bf R E F E R E N C E S}

\br\item{1. }R. V. Ambartzumian, ``Factorization
Calculus and Geometrical Probability", Cambridge
Univ. Press, Cambridge, 1990.
          
\item{2. }R. H. Aramyan, ``Curvature radii of
planar projections of convex bodies in $R^n$" [in
Russian], Izv. Akad. Nauk Armenii. Matematika,
[English translation: Journal of Contemporary
Math. Anal. (Armenian Academy of Sciences)],
vol. 37, no. 1, pp. 2 -- 14, 2002.

\item{3. }R. H. Aramyan, ``Generalized Radon
transform with an application in convexity theory"
[in Russian],
Izv. Akad. Nauk Armenii. Matematika, [English
translation: Journal of Contemporary Math. Anal.
(Armenian Academy of Sciences)], vol. 38, no. 3,
2003.

\item{4. }I. Ya. Bakelman, A. L. Verner, B. E. Kantor,
``Differential Geometry in the Large"[in Russian],
Nauka, Moskow, 1973.
       
\item{5. }W. Blaschke, ``Kreis und Kugel"
(Veit, Leipzig), 2nd Ed. De Gruyter, Berlin, 1956.

\item{6. } A. V. Pogorelov, ``Exterior Geometry of
Convex Surfaces" [in Russian], Nauka, Moscow, 1969.

\item{7. }A. P. Prudnikov, Yu. A. Brichov, O. N.
Marichev, ``Integrals and Series" [in Russian],
Nauka, Moscow, 1981.

\item{8. }W. Wiel, R. Schneider, ``Zonoids and
related Topics", in Convexity and its Applications,
Ed. P. Gruber and J. Wills, Birkhauser, Basel, 1983.

\hfill Institute of Mathematics

\hfill  Armenian Academy of Sciences

\hfill e.mail: rafik\@instmath.sci.am

\bye